\documentclass[10pt,a4paper,oneside]{amsart}

\usepackage{amsfonts,amssymb}
\usepackage{amsmath,amsthm}
\usepackage{graphicx, color}
\usepackage{enumerate,array,enumitem}
\usepackage{url}
\usepackage{tikz-cd}
\usepackage[hidelinks]{hyperref}

\usepackage[all]{xy}

\usetikzlibrary{arrows.meta,calc}

\newtheorem{thmA}{Theorem}

\newtheorem{theorem}{Theorem}[section]
\newtheorem{prop}[theorem]{Proposition}
\newtheorem{lemma}[theorem]{Lemma}
\newtheorem{cor}[theorem] {Corollary}

\newtheorem*{claim*}{Claim}

\theoremstyle{remark}
\newtheorem{remark}[theorem]{Remark}

\theoremstyle{definition}
\newtheorem{definition}[theorem]{Definition}

\newcommand{\id}{\textnormal{id}}

\newcommand{\Z}{\mathbb{Z}}

\newcommand{\C}{\mathbb{C}}
\newcommand{\R}{\mathbb{R}}

\newcommand{\metricdefect}[1][h]{\nu_{#1}}

\makeatletter
\newcommand{\xRightarrow}[2][]{\ext@arrow 0359\Rightarrowfill@{#1}{#2}}
\makeatother

\title[Discrete exterior calculus on surfaces]{Convergence of discrete exterior calculus for the Hodge--Dirac and Hodge--Laplace problems on surfaces}
\author{Thomas A. Wasserman}
\address{Department of Mathematics and Systems Analysis, Aalto University}
\email{thomas.wasserman@aalto.fi}
\date{\today}
\begin{document}
	\begin{abstract}
		We prove that discrete exterior calculus converges for the
		Hodge--Dirac and Hodge--Laplace problems on every closed oriented
		Riemannian surface. The proof applies the error bounds of	a companion paper on finite element exterior calculus. We prove and use that every closed oriented surface, with any metric, carries uniformly regular families of triangulations that are acute with a uniform
		margin.
	\end{abstract}
	\maketitle

\section{Introduction}
The Hodge--Dirac and Hodge--Laplace problems are variational formulations
of
\[
(d-{\star}d{\star})\alpha=\gamma\qquad\text{and}\qquad{-}(d{\star}d{\star}+{\star}d{\star}d)\alpha=\gamma
\]
where $\alpha$ and $\gamma$ are differential forms on a closed oriented Riemannian surface and $\star$ is its Hodge star. These variational formulations allow for any functional on the domain of $d$ to act as source. We show that discrete exterior calculus (DEC) converges for the Hodge--Dirac and Hodge--Laplace problems for any such source. We use the error bounds of \cite{WassermanFEEC} for finite element exterior
calculus approximations.  To apply them, we use that
on uniformly acute triangulations the DEC and finite element inner products are
uniformly equivalent, with a controlled difference \cite{ZhuChristiansenHuHirani2025}.
We show that every closed oriented surface, with any metric, carries
these acute triangulations.

Discrete exterior calculus on a triangulation approximates differential forms by simplicial
cochains, and the Hodge star by a diagonal matrix built geometrically \cite{Hirani2003,DesbrunHiraniLeokMarsden2005}.  The
inverse star is then diagonal too, so the codifferential and the Hodge Laplacian
are sparse matrices, which is a computational advantage over finite element
exterior calculus. On functions, the DEC Laplacian is the cotangent Laplacian, central to geometry processing \cite{MacNeal1949,Crane2013}. The diagonal star also makes the method local:
its Hodge star and codifferential at a simplex involve only that simplex and its
neighbours. This locality makes the method suited to discrete models in mathematical approaches to quantum physics. Those applications require having a convergence result for arbitrary sources.

Convergence of the DEC approximation has so far been proved for functions, and for
forms on flat spaces.  For functions on a surface, convergence for the Hodge--Laplace problem can be
established through linear finite elements, which share the cotangent Laplacian
with discrete exterior calculus \cite{Dziuk1988,Hildebrandt2006}.  In any
dimension, Schulz and Tsogtgerel proved convergence for functions directly, on
well-centred triangulations of polytopes \cite{SchulzTsogtgerel2020}.  For forms of
positive degree in two dimensions, Zhu, Christiansen, Hu and Hirani proved
convergence for the Hodge--Laplace problem on planar domains, on Delaunay
triangulations, by comparison with finite element exterior calculus
\cite{ZhuChristiansenHuHirani2025}.  In any dimension, convergence on well-centred
triangulations of polytopal domains has been proved for the Hodge--Laplace problem
\cite{GuzmanPotu2025} and, on topologically trivial domains, for the Hodge--Dirac
problem \cite{DabeticHiptmair2025}.

Convergence on a curved surface for arbitrary sources was out of reach for three reasons. First, discrete exterior calculus sees the metric of a
curved surface only through the edge lengths of the triangulation, so a geometric
error enters alongside the discretisation error; the analyses of forms
\cite{ZhuChristiansenHuHirani2025,GuzmanPotu2025,DabeticHiptmair2025} are set on
flat domains, where this error does not arise.  Second, their error bounds need more regularity of the source
than square integrability.  Third, the DEC inner product on cochains
differs from the continuum inner product by a variational crime that
does not tend to zero in the operator norm, so the theory of variational crimes
\cite{Holst2012} gives no convergence.  The companion paper \cite{WassermanFEEC} removes all three.  It
proves convergence of finite element exterior calculus on closed Riemannian
manifolds, with the piecewise flat metric of the triangulation in place of
the Riemannian metric and for every such source.  Its abstract error bound sees the
discrete inner products only through their deviation from the exact ones
relative to the norms of the forms and of their exterior derivatives.
Measured that way, the DEC inner product differs from the square-integrable
one by a term of the order of the mesh size on uniformly acute
triangulations, where Zhu, Christiansen, Hu and Hirani bound its difference
from the finite element inner product \cite{ZhuChristiansenHuHirani2025}.
We show that every surface carries such triangulations.

\subsection{Results}\label{sec: main results}
Throughout, $(\Sigma,g)$ is a closed oriented Riemannian surface.  We work with
smooth triangulations $K$ of $\Sigma$, whose closed simplices are smoothly
embedded, and with families of them that are \emph{uniformly regular} in the
sense of \cite{WassermanFEEC}: the embeddings of the simplices stay under
control across the family, and the mesh size $\ell_K$, the largest diameter
of a simplex of $K$, takes arbitrarily small values.  The metric enters only
through the lengths of the edges.  These lengths define the \emph{piecewise
flat metric} $g_K$, which makes each triangle a Euclidean triangle whose edges
have the lengths they have in $(\Sigma,g)$.

The Whitney map \cite{Whitney1957,Dodziuk1976} identifies the cochains of
discrete exterior calculus with Whitney forms, the lowest-order finite element
forms, so discrete exterior calculus is a finite element method with an inner
product of its own.  The \emph{DEC inner product} \cite{Hirani2003} of two
Whitney $k$-forms is the sum, over the $k$-simplices $\sigma$, of the product
of their integrals over $\sigma$, weighted by the ratio of the volume of the
circumcentric dual cell of $\sigma$ to the volume of $\sigma$, both measured
in $g_K$.  This needs the dual volumes to be positive, which holds when $K$ is
\emph{acute}, that is, when its triangles are acute in $g_K$. We use a form of \emph{pointwise acuteness}, checked at every point of the surface, because this pointwise condition survives subdivision and conformal changes of the metric. A family is called \emph{uniformly pointwise acute} if there is a margin of acuteness that is uniform in $K$.  The fine members of such a family are acute, with angles bounded below.

Our main theorem is that on such families discrete exterior calculus
converges.  We measure forms in the \emph{graph norm}
$(\|\alpha\|^2+\|d\alpha\|^2)^{1/2}$, and write $i_K$ for the inclusion of
the Whitney forms into the square-integrable forms.
\begin{thmA}[{Theorem~\ref{thm: dec convergence} and Corollary~\ref{cor: dec L2 convergence}}]\label{thm: intro b convergence}
	Let $(\Sigma,g)$ be a closed oriented Riemannian surface and $\mathcal K$ a uniformly regular
	and uniformly pointwise acute family of smooth triangulations of $\Sigma$, each carrying the Whitney forms with the DEC inner
	product of its piecewise flat metric.  Solve the Hodge--Dirac or the
	Hodge--Laplace problem with a source $F$, a bounded linear functional on the
	domain of the exterior derivative, and its discrete counterpart with source
	$F\circ i_K$.  Then the image under $i_K$ of the discrete solution converges to
	the solution in the graph norm as $\ell_K\to0$.  If $F$ is the inner product
	with a square-integrable form $\gamma$ whose codifferential $d^*\gamma$ is square
	integrable too, then the distance between them is at most of order
	$\ell_K(\|\gamma\|+\|d^*\gamma\|)$.
\end{thmA}
The proof applies the abstract error bound
of \cite{WassermanFEEC}.  The estimate of \cite{ZhuChristiansenHuHirani2025} for the
difference of the DEC and finite element inner products
makes their graph norms equivalent on the fine
members of such a family, so
discrete exterior calculus enters that bound only through its \emph{metric defect}, the
deviation of the DEC inner product from the square-integrable one relative to
the graph norm.  The agreement of the two inner products on closed Whitney
forms on each triangle makes the metric defect of the order of the mesh size.

The hypothesis of Theorem~\ref{thm: intro b convergence} can be met on every
surface.
\begin{thmA}[{Corollary~\ref{thm: well centred surfaces}}]\label{thm: intro b existence}
	Every closed oriented Riemannian surface carries a
	uniformly regular and uniformly pointwise acute family of smooth
	triangulations of arbitrarily small mesh size.
\end{thmA}
We employ Freudenthal subdivision, which from a parametrised base triangulation provides a uniformly regular family \cite{WassermanFEEC}. Since it copies the reference triangle in the coordinates of each base triangle, the subdivisions of a pointwise acute base triangulation form a uniformly pointwise acute family. Hence, we need to provide a pointwise acute base triangulation. Pointwise acuteness is unchanged by conformal changes of the metric, so on the sphere and the torus it suffices to find a base triangulation for a metric of constant curvature in the conformal class of $g$: we take the radially
projected icosahedron and an acute triangulation of the flat torus \cite{ItohYuan2009}.

Colin de~Verdi\`ere and Marin proved that every closed Riemannian surface has
acute geodesic triangulations of arbitrarily small mesh
\cite{ColinDeVerdiereMarin1990}. Their angles are measured at the corners,
between geodesic edges, a condition weaker than pointwise acuteness. In genus at least two we prove a parametrised analogue of their theorem, with pointwise acuteness in place of their optimal angle bounds. We use their construction: a metric close to $g$ that is conformal to a flat metric with cone points, in which the surface is cut into equilateral
triangles.

\label{sec: intro higher dimensions}In dimension three and above, the equivalence
of the two inner products extends by a scaling argument \cite{GuzmanPotu2025},
but the bound on the metric defect uses
the agreement of the two inner products on closed Whitney forms on each triangle,
so our argument stops at surfaces.  Nor do we know, on a curved manifold of
dimension three or more, of families of triangulations with the analogue of
uniform acuteness; the known results on acute and well-centred tetrahedral
meshes do not settle the question
\cite{VanderZee2013,KopczynskiPakPrzytycki2012}.

\subsection{Outline}
Section~\ref{sec: recollections} recalls what we need from
\cite{WassermanFEEC} and sets up discrete exterior calculus on a surface.
Section~\ref{sec: dec} defines pointwise acuteness, compares the
DEC and finite element inner products and bounds the metric defect, and proves
Theorem~\ref{thm: intro b convergence}.  Section~\ref{sec: acute triangulations}
provides pointwise acute base triangulations and proves
Theorem~\ref{thm: intro b existence}.

\section{Preliminaries}\label{sec: recollections}
The error estimates used here were established in a companion paper \cite{WassermanFEEC}. We recall the relevant results from it in Section~\ref{sec: recollections abstract}. We then set up discrete exterior calculus on a closed surface in Section~\ref{sec: dec complex}. Throughout, the scalars are real.

\subsection{Hilbert complexes and their approximations}\label{sec: recollections abstract}
We work in the setting of Hilbert complexes. We recall their definition, the formulation of the Hodge--Dirac and Hodge--Laplace problems on them, and their approximations.

\subsubsection{Hilbert complexes}
A \emph{Hilbert complex}
$(H^\bullet,V^\bullet,d)$ consists of Hilbert spaces $H^k$ and closed
densely defined operators $d\colon V^k\subseteq H^k\to H^{k+1}$ with
$d\circ d=0$ \cite{BruningLesch1992}; we write $\|\cdot\|$ for the norm of $H^\bullet$ and
$\|v\|_V^2=\|v\|^2+\|dv\|^2$ for the graph norm.  We write $Z^k=\ker d$,
$B^k=\operatorname{im}d$ and $\mathfrak H^k=Z^k\cap(B^k)^\perp$ for the
harmonic forms, and $V_\bullet$ for the domain of the adjoint
$d^*$ of $d$.  The complex is \emph{closed} when $B^k$ is closed in $H^k$
for every $k$, and then the Poincar\'e inequality
$\|v\|_V\le c_H\|dv\|$ holds on the orthogonal complement of $Z^\bullet$ in
$V^\bullet$, with a constant $c_H\ge1$.

\subsubsection{The Hodge--Dirac and Hodge--Laplace problems}

Both problems are posed in mixed weak form.  We write $(V^\bullet)^*$ for the dual
of $V^\bullet$ with its graph norm.  Given $F\in(V^\bullet)^*$, the
\emph{Hodge--Dirac problem} asks for
$(\alpha,\kappa)\in V^\bullet\times\mathfrak H^\bullet$ with
\begin{equation}\label{eq: Hodge Dirac problem}
	\begin{aligned}
		\langle \alpha, \lambda \rangle &= 0
		& \forall \lambda \in \mathfrak{H}^\bullet,\\
		\langle d \alpha, \beta\rangle+ \langle \alpha,d\beta\rangle
		+ \langle \kappa, \beta\rangle &= F(\beta)
		& \forall \beta\in V^{\bullet}.
	\end{aligned}
\end{equation}
For a source $\gamma\in H^\bullet$, that is,
$F=\langle\gamma,\cdot\rangle$, this is
\cite[(1)]{LeopardiStern2016}.  The problem is well-posed on a closed Hilbert
complex \cite[(11)]{WassermanFEEC}.  Given $F\in(V^k)^*$, the \emph{mixed
Hodge--Laplace problem} in degree $k$ asks for
$(\eta,\alpha,\kappa)\in V^{k-1}\times V^k\times\mathfrak H^k$ with
\begin{equation}\label{eq: Hodge Laplace problem}
    \begin{aligned}
    \langle \alpha, \lambda\rangle &=0 &\forall \lambda\in \mathfrak{H}^k,\\
    \langle \eta, \theta\rangle - \langle\alpha,d\theta\rangle&=0
        & \forall\theta\in V^{k-1},\\
    \langle d\eta, \beta\rangle +\langle d\alpha,d\beta\rangle
        +\langle \kappa, \beta\rangle& = F(\beta)
        & \forall \beta \in V^k .
    \end{aligned}
\end{equation}
For $F=\langle\gamma,\cdot\rangle$ with $\gamma\in H^k$ this is \cite[(18)]{arnold2010finite}. It is well-posed
\cite[(14)]{WassermanFEEC}.

\subsubsection{Approximations and error bounds}

The discrete complexes are compared with the continuous one through maps in
both directions.  An \emph{approximation} of a Hilbert complex
$(H^\bullet,V^\bullet,d)$ consists of finite-dimensional Hilbert complexes
$(H_h^\bullet,V_h^\bullet,d_h)$, indexed by $h$ in some set, each with inner
products $\langle\cdot,\cdot\rangle_h$ of its own, together with chain maps
$i_h\colon(V_h^\bullet,d_h)\to(V^\bullet,d)$ and
$\pi_h\colon(V^\bullet,d)\to(V_h^\bullet,d_h)$ of the domain complexes.
These satisfy $\pi_hi_h=\id$, and their operator norms for the graph norms
are bounded independently of $h$.  The map $i_h$
need not be isometric.  The discrete Hodge--Dirac and mixed
Hodge--Laplace problems are \eqref{eq: Hodge Dirac problem} and
\eqref{eq: Hodge Laplace problem} with $d$, $\langle\cdot,\cdot\rangle$ and
$\mathfrak H^\bullet$ replaced by $d_h$, $\langle\cdot,\cdot\rangle_h$ and the
harmonic forms $\mathfrak H_h^\bullet$ of the discrete complex, and with the
discrete source $F\circ i_h$.  The failure of $i_h$ to be isometric is measured by
the \emph{metric defect} in the graph norm, defined with the adjoint $i_h^*$ of
$i_h$ for the inner products of $H_h^\bullet$ and $H^\bullet$,
\begin{equation}\label{eq: graph metric defect}
	\metricdefect
	:=\sup_{0\ne x_h,y_h\in V_h^\bullet}
	\frac{|\langle(\mathrm{id}-i_h^*i_h)x_h,y_h\rangle_h|}
	{\|x_h\|_{V_h}\,\|y_h\|_{V_h}} .
\end{equation}

On an approximation of a closed Hilbert complex the discrete problems have unique solutions $(\alpha_h,\kappa_h)$ and $(\eta_h,\alpha_h,\kappa_h)$ \cite{WassermanFEEC}.  The
\emph{error} of the Hodge--Dirac problem is
$\|\alpha-i_h\alpha_h\|_V+\|\kappa-i_h\kappa_h\|$, and that of the
Hodge--Laplace problem is the sum of the three corresponding terms.
Regularity of the solution turns the
error bounds into a rate.  Let $(\mathcal W^\bullet,\|\cdot\|_{\mathcal W})$ be a
regularity scale for the approximation with scale parameter $\ell_h$
\cite[Definition~2.28]{WassermanFEEC} (a subspace of forms on which
$i_h\pi_h$ approximates the identity in $H^\bullet$ at the rate $\ell_h$) that satisfies a Gaffney estimate
\cite[Definition~2.30]{WassermanFEEC}, an abstract form of Gaffney's
inequality \cite[Corollary~2.1.6]{Schwarz1995}.  If $F=\langle\gamma,\cdot\rangle$ with
$\gamma$ in the domain $V_\bullet$ of $d^*$, then the error of either problem is at
most
\begin{equation}\label{eq: abstract rate recalled}
	C\big(\ell_h+\metricdefect\big)\big(\|\gamma\|+\|d^*\gamma\|\big),
\end{equation}
with $C$ depending only on the Poincar\'e constant, the bounds for $i_h$
and $\pi_h$ and the constants of the scale \cite[Theorem~2.35]{WassermanFEEC}.  For
a general source the error still tends to zero as $\ell_h+\metricdefect\to0$
\cite[Corollary~2.36]{WassermanFEEC}.  Of the terms of
\eqref{eq: abstract rate recalled} only the metric defect depends on the
discrete inner products.
\subsection{Discrete exterior calculus on a closed surface}\label{sec: dec complex}

From here on $(\Sigma,g)$ is a closed oriented Riemannian surface.  We
recall the $L^2$ de~Rham complex of $\Sigma$ and its triangulations, and then set up discrete exterior calculus on a triangulation.

\subsubsection{The $L^2$ de~Rham complex}

The Hodge star and integration give the inner
product $\langle\alpha,\beta\rangle=\int_\Sigma\alpha\wedge\star\beta$ on
smooth $k$-forms; $L^2\Lambda^k(\Sigma)$ is its completion,
$V^\bullet$ the domain of the weak exterior derivative
\cite[(2.3) and (2.6)]{Licht2023smoothed} and
$V_\bullet$ the domain of its adjoint $d^*$.  The complex
$(L^2\Lambda^\bullet(\Sigma),V^\bullet,d)$, the \emph{$L^2$ de~Rham
complex}, is a closed Hilbert complex \cite[Section~7]{Licht2023smoothed}.  Rates need the solution to have one
derivative, that is, to lie in the Sobolev space $H^1\Lambda^k(\Sigma)$ of forms with
square-integrable first covariant derivative, with norm
$\|\omega\|_{H^1}=\|\omega\|+\|\nabla\omega\|$
\cite[(2.7)]{Licht2023smoothed}.

\subsubsection{Triangulations and piecewise flat metrics}
A \emph{smooth triangulation} $K$ of $\Sigma$ is a collection of smooth
simplices, each given by a smooth embedding $\imath_T\colon\Delta^k\to \Sigma$ of
the reference $k$-simplex, compatible along subsimplices
\cite[Section~3.1]{Licht2023smoothed}.  We write $h_T$ for the $g$-diameter of
a simplex $T$ of $K$ and $\ell_K:=\max_Th_T$ for the largest of these.  A
family $\mathcal K$ with $\inf_{K\in\mathcal K}\ell_K=0$ is \emph{uniformly
regular} when the constants of regularity of $K$ (see \cite[Definition~3.2]{WassermanFEEC} and \cite[Section~3.2]{Licht2023smoothed}) are bounded independently of $K$.

The \emph{piecewise flat metric} $g_K$ of $K$ is the metric on the
interiors of the triangles that is flat on each of them and gives each
edge its $g$-length.  It exists when, for every triangle, these
lengths are the edge lengths of a non-degenerate Euclidean triangle; we then
call $K$ \emph{realisable} \cite[Definition~3.3]{WassermanFEEC}.

From now on, we let $K$ be a
realisable smooth triangulation of $\Sigma$ with piecewise flat metric $g_K$.
Every triangle $T$ of $K$ with the metric $g_K$ is isometric to a
Euclidean triangle, by an affine map composed with $\imath_T^{-1}$, and we
identify the two.  All lengths, areas and circumcentres below are taken in
these Euclidean triangles.  An edge lying in two triangles has the same
length in both, so the volume $|\sigma|$ of a simplex $\sigma$, its length
for an edge and its area for a triangle, does not depend on the choice; the
volume of a vertex is $1$. The orientation of $\Sigma$ orients the triangles, and we fix an orientation of every edge.

\subsubsection{Dual cells}\label{sec: dec dual cells}
In discrete exterior calculus the metric $g_K$ enters only through the discrete
Hodge star, which is built using
circumcentric dual cells.  From here on every
triangle $T$ of $K$ is acute, so that its circumcentre $c(T)$, the centre of
its circumcircle, lies in its interior \cite[Section~4]{VanderZee2010}. The \emph{dual cell} of a simplex is a cell of
complementary dimension built from the circumcentres of the simplex and of the
simplices containing it \cite[Definition~2.4.5]{Hirani2003}.  The
dual cell $\star T$ of a triangle $T$ is the point $c(T)$.  The dual cell
$\star\sigma$ of an edge or a vertex $\sigma$ is the union of its parts
$\star\sigma\cap T$ in the triangles $T$ containing $\sigma$.  If $\sigma$ is
an edge of $T$, the part $\star\sigma\cap T$ is the segment from the midpoint
of $\sigma$ to $c(T)$.  If $\sigma$ is a vertex of $T$, it is the quadrilateral whose vertices are
$\sigma$, the midpoints of the two edges of $T$ at $\sigma$, and $c(T)$; the
quadrilaterals of the three vertices of $T$ tile $T$.  We write $|\star\sigma\cap T|$ for
the volume of $\star\sigma\cap T$, so that $\star\sigma$ has volume
$|\star\sigma|=\sum_{T\supset\sigma}|\star\sigma\cap T|$, and for a
triangle $T$ we set $|\star T\cap T|=|\star T|:=1$.  Every $|\star\sigma\cap T|$, and hence every $|\star\sigma|$, is nonzero, and
so positive, because $c(T)$ lies in the interior of $T$.

\subsubsection{Cochains and the DEC inner product}
Discrete exterior calculus, in the form given to it by Hirani, takes the simplicial cochain complex
$(C^\bullet(K),d)$ for its discrete forms and exterior derivative
\cite[Definitions~3.2.3 and~3.6.5]{Hirani2003}, and adds a discrete Hodge star
\cite[Definition~4.1.1]{Hirani2003}.

The Hodge star of a $k$-form is a $(2-k)$-form, so the discrete Hodge star
sends a $k$-cochain $a$ to a function $\star a$ on the dual cells of
dimension $2-k$.  The inner product then pairs a cochain with the
Hodge star of another, as $\int_\Sigma\alpha\wedge\star\beta$ does.

\begin{definition}[{\cite[Definitions~3.9 and~3.12]{GuzmanPotu2025}}]\label{def: dec complex inner product}
    The \emph{discrete Hodge star} of $a\in C^k(K)$ is the function
    $\star a$ on the circumcentric duals of the $k$-simplices with
    \[
        (\star a)(\star\sigma)=\frac{|\star\sigma|}{|\sigma|}\,a(\sigma).
    \]
    The \emph{DEC inner product} of
    $a,b\in C^k(K)$ is
    \begin{equation}\label{eq: dec complex inner product}
        \langle a,b\rangle_D:=\sum_\sigma a(\sigma)\,(\star b)(\star\sigma)
        =\sum_\sigma\frac{|\star\sigma|}{|\sigma|}\,a(\sigma)\,b(\sigma),
    \end{equation}
    the sum over the $k$-simplices of $K$.
\end{definition}

In the bases given by the simplices and by their duals, the Hodge star is
the diagonal matrix with entries $|\star\sigma|/|\sigma|$.  The form $\langle\cdot,\cdot\rangle_D$ is symmetric and
bilinear, and it is an inner product because every $|\star\sigma|$ is
positive.

\subsubsection{Whitney forms}

We carry discrete exterior calculus over to forms, where it can be compared
with finite elements, through the Whitney map.

\begin{lemma}[{\cite[Definition~1.1, (1.4) and (1.7)]{Dodziuk1976}, \cite[Section~4.3]{Arnold2006}}]\label{lem: whitney map}
    Integration over the simplices, which sends $\alpha\in V_K^k$ to the
    cochain $\sigma\mapsto\alpha(\sigma):=\int_\sigma\alpha$, is an isomorphism of
    cochain complexes from the lowest-order finite element complex
    $V_K^\bullet:=\mathcal P_1^-\Lambda^\bullet(K)$
    \cite[Definition~3.18]{WassermanFEEC} onto $C^\bullet(K)$.
\end{lemma}
The inverse to this isomorphism is called the \emph{Whitney map}. It carries the DEC inner product over to forms: for
$\alpha,\beta\in V_K^k$,
\begin{equation}\label{eq: dec inner product forms}
    \langle\alpha,\beta\rangle_D
    :=\sum_\sigma\frac{|\star\sigma|}{|\sigma|}\,\alpha(\sigma)\,\beta(\sigma),
\end{equation}
the sum over the $k$-simplices of $K$.  We write $H_{K,D}^k$ for $V_K^k$
with this inner product, $\|\cdot\|_D$ for its norm and $\|\cdot\|_{V_{K,D}}$
for its graph norm.  Discrete exterior calculus is then the Hilbert complex
$(H_{K,D}^\bullet,V_K^\bullet,d)$.

In contrast, on $V_K^\bullet$ the \emph{finite element inner product} of $g_K$ is
$\langle\alpha,\beta\rangle_K=\int_\Sigma\alpha\wedge\star_K\beta$, where
$\star_K$ is the Hodge star of $g_K$.  We write $H_K^\bullet$ for
$V_K^\bullet$ with this inner product, $\|\cdot\|_K$ for its norm and
$\|\cdot\|_{V_K}$ for its graph norm.
The results of \cite{WassermanFEEC} apply to the finite element complex
$(H_K^\bullet,V_K^\bullet,d)$, but not directly to $(H_{K,D}^\bullet,V_K^\bullet,d)$, because $\langle\cdot,\cdot\rangle_D$ is not the $L^2$
inner product of any metric on the triangles.\footnote{Indeed, let $v\ne v'$ be
vertices of a triangle $T$, and let $\lambda_v,\lambda_{v'}\in V_K^0$ be
their hat functions, which on $T$ are barycentric coordinates
\cite[(4.3)]{Arnold2006}.  Each vanishes at every vertex but its own, so
$\langle\lambda_v,\lambda_{v'}\rangle_D=0$ by
\eqref{eq: dec inner product forms}.  But $\int_T\lambda_v\lambda_{v'}>0$
for every metric.} To apply \cite{WassermanFEEC}, we compare the DEC and finite element inner products.

\section{Convergence of discrete exterior calculus}\label{sec: dec}

This section compares the DEC inner product with the finite element one and
proves that discrete exterior calculus converges.  The comparison
needs acute triangulations whose angles are bounded below uniformly, and every
surface carries a uniformly regular family of this kind
(Corollary~\ref{thm: well centred surfaces}).  On these triangulations the
difference of the two inner products is bounded in the graph norm by the mesh
size \cite{ZhuChristiansenHuHirani2025}.
Section~\ref{sec: dec inner products surface} recalls this bound and
uses it to show that, on fine meshes, the two graph norms are equivalent and
the metric defect is of the order of the mesh size.
Section~\ref{sec: dec convergence} then proves the convergence theorem.

\subsection{Pointwise acute families}\label{sec: dec well centred families}
In Section~\ref{sec: dec inner products surface} we need a family of
triangulations that are acute with angles away from $0$.  This
subsection gives a pointwise condition on smooth triangulations that guarantees
this, and subdivides a base triangulation that satisfies it into a uniformly pointwise acute family
(Proposition~\ref{prop: well centred from base}).

\subsubsection{Pointwise acuteness}\label{sec: dec norm equivalence}\label{sec: dec acute base angles}

A realisable smooth triangulation $K$ is \emph{acute} if all its triangles are
acute in $g_K$.  For an infinite family we need a notion of uniform
acuteness.  Zhu, Christiansen, Hu and Hirani call a family
\emph{uniformly acute} if one upper bound less than $\pi/2$ holds for all
angles of all its triangles
\cite[Section~2.1, condition~(ii)]{ZhuChristiansenHuHirani2025}.  This also
bounds all angles below uniformly
\cite[Remark~3.1]{ZhuChristiansenHuHirani2025}, which makes the comparison of
the two inner products in Section~\ref{sec: dec inner products surface} uniform
over the family.  We impose such a
bound at every point of every triangle instead, because the pointwise condition
passes to Freudenthal subdivisions
(Proposition~\ref{prop: well centred from base}) and is unchanged by conformal
changes of $g$ (Section~\ref{sec: dec well centred surfaces}).  For a triangle
$T$ of a smooth triangulation $K$ and $x\in\Delta^2$, we write
$(\imath_T^*g)_x$ for the value at $x$ of the pullback metric $\imath_T^*g$ on
the reference triangle $\Delta^2\subset\R^2$. For each such $x$ we regard $\Delta^2$ as the convex hull of the origin and the basis vectors in the Euclidean plane $(\R^2,(\imath_T^*g)_x)$, and measure its angles there.

\begin{definition}\label{def: dec uniformly acute}\label{def: acute base angles}
    A smooth triangulation $K$ of $\Sigma$ is \emph{pointwise acute} if, for
    every triangle $T$ of $K$ and every $x\in\Delta^2$, the triangle $\Delta^2$ is
    acute in $(\imath_T^*g)_x$.  A family $\mathcal K$ of smooth triangulations
    of $\Sigma$ is \emph{uniformly pointwise acute}, with margin $\delta>0$, if
    for every $K\in\mathcal K$, every triangle $T$ of $K$ and every
    $x\in\Delta^2$, every angle of $\Delta^2$ in $(\imath_T^*g)_x$ is at most
    $\pi/2-\delta$.
\end{definition}

Acuteness is a property of the flat triangles of $g_K$, pointwise acuteness one
of the parametrised triangles in $g$.  On fine meshes the second implies the
first (Lemma~\ref{lem: dec acute flat} below), because nearby inner products give
nearby angles.  Using \cite[Theorem~2.1]{LinSinnamon2012} one shows that this
holds in the following form.

\begin{lemma}\label{lem: angle comparison}
    Let $G$ and $G'$ be inner products on a real vector space, and let
    $0<a\le b$ with $a^2G\le G'\le b^2G$.  Then for all nonzero $u$ and $v$
    the angle between $u$ and $v$ measured in $G'$ differs from the angle
    measured in $G$ by at most
    \begin{equation}\label{eq: angle distortion}
        2\arcsin\frac{b-a}{b+a}.
    \end{equation}
\end{lemma}

Comparing the two notions of acuteness uses the metric error $\epsilon_K$
of $K$, the smallest $\epsilon$ with $(1-\epsilon)g\le g_K\le(1+\epsilon)g$
almost everywhere \cite[Definition~3.3]{WassermanFEEC}.

\begin{lemma}\label{lem: dec acute flat}
    Let $K$ be a realisable smooth triangulation of $\Sigma$ that is pointwise
    acute with margin $\delta>0$.  If $\epsilon_K$ is small enough, depending
    only on $\delta$, then $K$ is acute, with all angles in $g_K$ at least
    $\delta$.
\end{lemma}

\begin{proof}
Let $T$ be a triangle of $K$, let $G_x:=(\imath_T^*g)_x$ for $x\in\Delta^2$, and
let $G':=\imath_T^*g_K$.  By the definition of
$\epsilon_K$, $(1-\epsilon_K)G_x\le G'\le(1+\epsilon_K)G_x$ for almost every
$x\in\Delta^2$, hence for every $x$, as $G_x$ depends continuously on $x$.
So, for $\epsilon_K<1$, Lemma~\ref{lem: angle comparison} applies: the angles of
$\Delta^2$ in $G'$ exceed those in $G_x$ by at most the
bound~\eqref{eq: angle distortion}, which
tends to $0$ as $\epsilon_K\to0$.  For $\epsilon_K$ small enough that bound is at most $\delta/2$, and the angles of $T$ in
$g_K$ are then at most $\pi/2-\delta/2$, hence also at least $\delta$.
\end{proof}
\subsubsection{Freudenthal subdivisions}

We obtain our families as Freudenthal subdivisions
\cite[Section~3.2.2]{WassermanFEEC}, which on a
surface copy the shapes of triangles. Let $\sigma$ be a
triangle in the plane.  Every triangle of the $j$-th Freudenthal subdivision of
$\sigma$ is a translate of $2^{-j}$ times a conjugate of $\sigma$
\cite[Section~3.2.2]{WassermanFEEC}.  The conjugates of $\sigma$ are
$\sigma$ itself and its image under a point reflection. So every triangle of the $j$-th Freudenthal subdivision
of $\sigma$ has the angles of $\sigma$, measured in any inner product on the plane;
in particular, if $\sigma$ is acute, then so are all of these
triangles.

\subsubsection{Pointwise acute families from a base triangulation}\label{sec: acute families on every surface}
Any smooth triangulation $K_0$ comes with a homeomorphism $\chi\colon |K_0|\to \Sigma$. This splits for each simplex $\sigma$ of $K_0$ its parametrisation as $\imath_\sigma=\chi \circ a_\sigma$, where $a_\sigma$ is affine \cite[Section~3.2.2]{WassermanFEEC}.

\begin{prop}\label{prop: well centred from base}
    Let $K_0$ be a pointwise acute smooth triangulation of $\Sigma$.  Then
    its Freudenthal subdivisions $K_j$, $j\ge0$,
    form a uniformly pointwise acute family.
\end{prop}

\begin{proof}
The angles of $\Delta^2$ in $(\imath_\sigma^*g)_y$ depend continuously on
$y\in\Delta^2$, which is compact, and $K_0$ has finitely many triangles
$\sigma$, so there is $\delta>0$ such that all of these angles are at most
$\pi/2-\delta$.  Let $\tau$ be a triangle of $K_j$, contained in the triangle
$\sigma$ of $K_0$, and let $T:=a_\sigma^{-1}(\tau)$, a triangle of the $j$-th
Freudenthal subdivision of $\Delta^2$.  Then $\imath_\tau=\imath_\sigma\circ A$ for
the affine map $A:=a_\sigma^{-1}\circ a_\tau$ of $\Delta^2$ onto $T$, so for
$x\in\Delta^2$ the angles of $\Delta^2$ in $(\imath_\tau^*g)_x$ are those of $T$ in
$(\imath_\sigma^*g)_{A(x)}$.  By the description of the Freudenthal subdivisions
above, they are the angles of $\Delta^2$ in $(\imath_\sigma^*g)_{A(x)}$, which are
at most $\pi/2-\delta$.  So the family is uniformly pointwise acute with
margin $\delta$.
\end{proof}

\subsection{The DEC and finite element inner products}\label{sec: dec inner products surface}
In this subsection we compare the DEC inner product with the finite element
inner product of $g_K$: their difference is bounded in the graph norm by the
mesh size \cite[Theorem~5.4]{ZhuChristiansenHuHirani2025}.
The bound plays two roles.  On fine meshes it makes the two graph norms
equivalent, which brings discrete
exterior calculus under the error bounds of
Section~\ref{sec: recollections abstract}.
Together with the corresponding bound for finite elements, it also controls
the metric defect, the only term of those bounds that depends on the inner
products.

\subsubsection{The difference of the two inner products}
We use the following bound on the difference of the
two inner products in the graph norm by the mesh size, here the largest
diameter $\ell_K$.  The bound is due to Zhu, Christiansen, Hu and
Hirani, who also survey the history of such comparisons
\cite[Sections~2.4 and~2.5]{ZhuChristiansenHuHirani2025}.

\begin{prop}[{\cite[Theorem~5.4]{ZhuChristiansenHuHirani2025}}]\label{prop: dec feec consistency}
    Let $\theta_0>0$ and suppose that $K$ is acute with angles at least $\theta_0$.  There is a constant $C$, depending
    only on $\theta_0$ and the shape constant of $K$
    \cite[Definition~3.2]{WassermanFEEC}, such that for all
    $\alpha,\beta\in V_K^\bullet$
    \[
        \big|\langle\alpha,\beta\rangle_D-\langle\alpha,\beta\rangle_K\big|
        \le C\,\ell_K\,\|\alpha\|_{V_K}\|\beta\|_{V_K}.
    \]
\end{prop}

\begin{proof}
This is \cite[Theorem~5.4]{ZhuChristiansenHuHirani2025}, which is stated for a
surface which is a piecewise linear surface in $\R^3$ or a planar domain. The surface
$(\Sigma,g_K)$ need not embed in $\R^3$ with flat triangles.  However, their proof sums over the triangles the bounds of their Lemmas~5.2 and~5.3, which only use the difference of the two inner products on a single triangle, so carries over.  Their mesh size, the length of the longest edge of $K$ in $g_K$, is at most a multiple of $\ell_K$ depending only on the shape
constant of $K$, because an edge of a triangle $T$ has $g_K$-length its
$g$-length, at most a fixed multiple of $\sup_{\Delta^2}|d\imath_T|$.
\end{proof}

On fine meshes the bound makes the graph norms of the two inner products
equivalent.  Let $K$, $\theta_0$ and $C$ be as in
Proposition~\ref{prop: dec feec consistency}, with $C\ell_K\le1/4$, and let
$\alpha\in V_K^\bullet$.  The squared graph norms differ by
\[
    \|\alpha\|_{V_{K,D}}^2-\|\alpha\|_{V_K}^2
    =\big(\|\alpha\|_D^2-\|\alpha\|_K^2\big)
     +\big(\|d\alpha\|_D^2-\|d\alpha\|_K^2\big).
\]
Proposition~\ref{prop: dec feec consistency} bounds the first term in
absolute value by $C\ell_K\|\alpha\|_{V_K}^2$, and the second by
$C\ell_K\|d\alpha\|_{V_K}^2$, where
$\|d\alpha\|_{V_K}=\|d\alpha\|_K\le\|\alpha\|_{V_K}$ because $dd\alpha=0$.
So the difference is at most
$2C\ell_K\|\alpha\|_{V_K}^2\le\frac12\|\alpha\|_{V_K}^2$ in absolute
value.  Hence
\begin{equation}\label{eq: dec graph norm equivalence}
    \tfrac12\|\alpha\|_{V_K}^2\le\|\alpha\|_{V_{K,D}}^2\le2\|\alpha\|_{V_K}^2
    \qquad\text{for }\alpha\in V_K^\bullet.
\end{equation}

\subsubsection{The metric defect}\label{sec: dec metric defect}

The metric defect compares the DEC inner product with the $L^2$ inner
product of $g$. The distance between the DEC inner product and the finite element inner product is bounded by Proposition~\ref{prop: dec feec consistency}, and between that and the $L^2$
inner product by \cite[Theorem~3.21]{WassermanFEEC}.

\begin{prop}\label{prop: dec graph defect}
    Let $\mathcal K$ be uniformly regular and uniformly pointwise acute with
    margin $\delta$, and give each $V_K^\bullet$ the DEC inner product.  Then
    for every $K\in\mathcal K$ with $\ell_K$ small enough, the metric defect
    \eqref{eq: graph metric defect} of the inclusion $i_K$ satisfies
    \[
        \metricdefect[K]\le C\,\ell_K,
    \]
    with $C$ and the threshold for $\ell_K$ depending only on
    $(\Sigma,g)$, $\delta$ and bounds for the constants of regularity of the
    family.
\end{prop}

\begin{proof}
Let $K\in\mathcal K$ with $\ell_K$ small enough that $K$ is
realisable \cite[Proposition~3.5]{WassermanFEEC} and acute with angles at least $\delta$
(Lemma~\ref{lem: dec acute flat}), that \cite[Theorem~3.21]{WassermanFEEC}
applies and that the norm equivalence \eqref{eq: dec graph norm equivalence}
holds.
Let $x,y\in V_K^\bullet$, and take the adjoint $i_K^*$ with respect to
$\langle\cdot,\cdot\rangle_D$.  Then
$\langle(\mathrm{id}-i_K^*i_K)x,y\rangle_D=\langle x,y\rangle_D-\langle i_Kx,i_Ky\rangle$,
so it suffices to bound the right-hand side by
$C\ell_K\|x\|_{V_{K,D}}\|y\|_{V_{K,D}}$.  To do so we compare both inner products
with the finite element inner product $\langle\cdot,\cdot\rangle_K$ of $g_K$:
\[
    \big|\langle x,y\rangle_D-\langle i_Kx,i_Ky\rangle\big|
    \le\big|\langle x,y\rangle_D-\langle x,y\rangle_K\big|
      +\big|\langle x,y\rangle_K-\langle i_Kx,i_Ky\rangle\big| .
\]
Proposition~\ref{prop: dec feec consistency}, with $\theta_0=\delta$, bounds
the first difference by $C\ell_K\|x\|_{V_K}\|y\|_{V_K}$, with $C$ independent
of $K$ by uniform regularity.

The second difference is $\langle(\mathrm{id}-i_K^{*_K}i_K)x,y\rangle_K$,
where $i_K^{*_K}$ is the adjoint taken with respect to
$\langle\cdot,\cdot\rangle_K$.
\cite[Theorem~3.21]{WassermanFEEC} bounds the operator norm of
$\mathrm{id}-i_K^{*_K}i_K$ on $H_K^\bullet$ by $C\ell_K$
\cite[(49)]{WassermanFEEC}.  So the second difference is at most
$C\ell_K\|x\|_K\|y\|_K\le C\ell_K\|x\|_{V_K}\|y\|_{V_K}$.
Equation \eqref{eq: dec graph norm equivalence}
converts the graph norms $\|\cdot\|_{V_K}$ into those of $\langle\cdot,\cdot\rangle_D$, at
the cost of a factor $2$.  Taking the supremum over $x$ and $y$ of unit DEC graph norm then gives the
assertion.
\end{proof}

\begin{remark}\label{rem: dec metric defect}
    The bound of Proposition~\ref{prop: dec graph defect} is in the graph
    norm, and the operator norm of $\mathrm{id}-i_K^*i_K$ on
    $H_{K,D}^\bullet$ need not tend to zero.  Let $K$ be a triangulation of a
    flat torus by equilateral triangles, so that $g_K=g$, whose vertices
    admit a proper $3$-colouring; such triangulations exist with arbitrarily
    small triangles.  Let $\alpha\in V_K^0$ take the values
    $1,-1,0$ on the three colour classes.  Then one
    checks that $\langle\alpha,\alpha\rangle_D=4\langle\alpha,\alpha\rangle$.  So
    $\langle(\mathrm{id}-i_K^*i_K)\alpha,\alpha\rangle_D
    =\langle\alpha,\alpha\rangle_D-\langle i_K\alpha,i_K\alpha\rangle
    =\tfrac34\langle\alpha,\alpha\rangle_D$, and
    $\|\mathrm{id}-i_K^*i_K\|\ge3/4$ however fine $K$ is.  An error bound
    containing $\|\mathrm{id}-i_K^*i_K\|\,\|\gamma\|$, such as
    \cite[Corollary~3.11]{Holst2012}, therefore does not show that discrete
    exterior calculus converges.
\end{remark}

\subsection{The convergence theorem}\label{sec: dec convergence}
This subsection shows that discrete exterior calculus converges for the
Hodge--Dirac and mixed Hodge--Laplace problems.  It combines the abstract error bounds of
Section~\ref{sec: recollections abstract} with the norm equivalence
\eqref{eq: dec graph norm equivalence} and the bound on the metric defect
(Proposition~\ref{prop: dec graph defect}).

\subsubsection{The rate of convergence}

We now prove the rate at which the DEC solutions converge to the
continuous solution for sources in the domain of the codifferential.

\begin{theorem}\label{thm: dec convergence}\label{thm: dec convergence hodge dirac}\label{thm: dec convergence hodge laplace}
    Let $\mathcal K$ be a uniformly regular family of smooth triangulations
    of the closed oriented surface $\Sigma$ that is uniformly pointwise acute with margin
    $\delta$, and give each
    $V_K^\bullet$ the DEC inner product.  Consider either the Hodge--Dirac
    problem \eqref{eq: Hodge Dirac problem} with source $F=\langle\gamma,\cdot\rangle$
    for some $\gamma\in V_\bullet$, or the mixed Hodge--Laplace problem
    \eqref{eq: Hodge Laplace problem} in degree $k$ with source
    $F=\langle\gamma,\cdot\rangle$ for some $\gamma\in V_k$.  Then for every
    $K\in\mathcal K$ with $\ell_K$ small enough (as in
    Proposition~\ref{prop: dec graph defect}), the discrete problem
    on $(H_{K,D}^\bullet,V_K^\bullet,d)$ with source $F\circ i_K$ has a unique
    solution, and its error $\mathrm{err}_K$
    (Section~\ref{sec: recollections abstract}) satisfies
    \begin{equation}\label{eq: dec rate}
        \mathrm{err}_K
        \le C\ell_K\big(\|\gamma\|+\|d^*\gamma\|\big).
    \end{equation}
    The constant $C$ and the threshold for $\ell_K$ depend only on
    $(\Sigma,g)$, $\delta$ and bounds for the constants of regularity of the
    family, and not on $K$ or on $\gamma$.
\end{theorem}

\begin{proof}
We apply the abstract error bound \cite[Theorem~2.35]{WassermanFEEC} to the
complexes $(H_{K,D}^\bullet,V_K^\bullet,d)$ with the inclusions $i_K$ and the
smoothed projections $\pi_K$ of \cite[Theorem~6.11]{Licht2023smoothed}.  For
the finite element inner products of $g_K$ its hypotheses hold by
\cite[Theorem~3.21]{WassermanFEEC},
\cite[Corollary~3.24]{WassermanFEEC} and
\cite[Lemma~3.16]{WassermanFEEC}, with constants depending only on
$(\Sigma,g)$ and bounds for the constants of regularity of the family.  These hypotheses
involve the discrete inner product only through the graph norms of
$V_K^\bullet$, so by the norm equivalence \eqref{eq: dec graph norm equivalence}
they hold for the DEC inner products as well.  So the discrete
problems have unique solutions, and
$\mathrm{err}_K\le C\big(\ell_K+\metricdefect[K]\big)\big(\|\gamma\|+\|d^*\gamma\|\big)$
by \eqref{eq: abstract rate recalled}. Proposition~\ref{prop: dec graph defect} then gives \eqref{eq: dec rate}.
\end{proof}
A family of triangulations satisfying the hypotheses of
Theorem~\ref{thm: dec convergence} can be found on every closed oriented
surface, see Corollary~\ref{thm: well centred surfaces} below.

\subsubsection{Sources in the dual of the domain}

For a general source in the dual of the domain the theorem gives no rate,
but discrete exterior calculus still converges, because the metric defect
tends to zero with the mesh size.

\begin{cor}\label{cor: dec L2 convergence}
    Let $\mathcal K$, $\Sigma$ and the DEC inner products be as in
    Theorem~\ref{thm: dec convergence}.  Consider either the Hodge--Dirac
    problem with source $F\in(V^\bullet)^*$, or the mixed
    Hodge--Laplace problem in degree $k$ with source $F\in(V^k)^*$, and for each $K$ the discrete source $F\circ i_K$.  Then for every $\varepsilon>0$ there is $\ell>0$ such
    that the error $\mathrm{err}_K$ of Theorem~\ref{thm: dec convergence}
    satisfies $\mathrm{err}_K\le\varepsilon$ for all $K\in\mathcal K$ with
    $\ell_K\le\ell$.
\end{cor}

\begin{proof}
This follows from \cite[Corollary~2.36]{WassermanFEEC} and
Proposition~\ref{prop: dec graph defect}.
\end{proof}

\section{Pointwise acute base triangulations of surfaces}\label{sec: acute triangulations}
Proposition~\ref{prop: well centred from base} and
\cite[Theorem~3.14]{WassermanFEEC} turn a pointwise
acute base triangulation (Definition~\ref{def: dec uniformly acute}) into a
uniformly regular, uniformly pointwise acute family.  In this section
we prove that the surface $(\Sigma,g)$ has such a base triangulation.

\begin{theorem}\label{prop: acute base exists}
    Every closed smooth oriented Riemannian surface has a pointwise acute smooth triangulation.
\end{theorem}

\begin{proof}[Proof of Theorem~\ref{prop: acute base exists}]
	It suffices to treat connected surfaces.  The rest of the section gathers the ingredients.  On the sphere and the
	torus we construct the base triangulation on the constant curvature model
	and carry it over to $(\Sigma,g)$ by a conformal diffeomorphism, which
	keeps pointwise acuteness
	(Section~\ref{sec: dec well centred surfaces}).
	In genus at least two the theorem is proved in
	Section~\ref{sec: cdvm appendix}.
\end{proof}

This provides the base triangulation to apply Proposition~\ref{prop: well centred from base} to.

\begin{cor}\label{thm: well centred surfaces}
    Let $(\Sigma,g)$ be a closed oriented Riemannian surface.  There is a
    pointwise acute smooth triangulation $(K_0,\chi)$ of $\Sigma$ whose
    Freudenthal subdivisions $K_j$, $j\ge0$, form a uniformly regular family of
    smooth triangulations of $\Sigma$ that is uniformly pointwise acute
    (Definition~\ref{def: dec uniformly acute}), with $\ell_{K_j}\to0$.
\end{cor}

\begin{proof}
By
Theorem~\ref{prop: acute base exists}, $\Sigma$ has a pointwise acute smooth
triangulation $(K_0,\chi)$, and by Proposition~\ref{prop: well centred from base}
its Freudenthal subdivisions form a uniformly pointwise acute family.  By
\cite[Theorem~3.14]{WassermanFEEC} they form a
uniformly regular family of smooth triangulations of $\Sigma$ with
$\ell_{K_j}\to0$.
\end{proof}

\begin{remark}
	Colin de Verdi\`ere and Marin prove that every closed Riemannian surface has
	acute \emph{geodesic} triangulations \cite{ColinDeVerdiereMarin1990}. They bound the angles at the corners of the triangles, between geodesic edges. Pointwise acuteness asks for more.  Each triangle $\sigma$ carries a parametrisation $\imath_\sigma$, and $\Delta^2$ must be
	acute in the inner product $(\imath_\sigma^*g)_x$ at every point $x$, not only at
	the corners.
	Finding acute triangulations is far from easy \cite[Remark~5.2]{GasterLoustauMonsaingeon2023}.  The existence proof of Colin de~Verdi\`ere and Marin is highly
	non-constructive. Brunck \cite{Brunck2022} constructs them for polygonal complexes
	of constant curvature.
\end{remark}

\subsection{The sphere and the torus}\label{sec: dec well centred surfaces}
We construct the base triangulation on a model of constant curvature.  Pointwise
acuteness is conformally invariant, and by uniformisation $(\Sigma,g)$ is
conformally equivalent to a surface of constant curvature: the round sphere in
genus zero and a flat torus $\C/\Lambda$ in genus one
\cite[Theorem~22 and Corollaries~3 and~4]{Donaldson2011}.

\subsubsection{The round sphere}\label{sec: acute sphere}
We show that the inverse of the radial projection of the unit sphere onto an
inscribed icosahedron makes the icosahedron a pointwise acute smooth
triangulation of the sphere.  It is clearly a smooth triangulation, so we need to check that each face is acute, at every one of its points, in the round
metric pulled back by the inverse projection.
This is the following computation.

Let $K_0$ be the icosahedron inscribed in the unit sphere with vertices the cyclic
permutations of $(0,\pm1,\pm\phi)/(\phi+2)^{1/2}$, where $\phi=\frac12(1+\sqrt5)$,
and let $\chi(p):=p/|p|$ be the radial projection onto the unit sphere.  On a face
$\sigma$ with centre $c$, $\chi$ is the dilation by $1/|c|$, which takes $\sigma$
into the plane touching the unit sphere at $c/|c|$, followed by the inverse of
the gnomonic projection of the unit sphere onto that plane
\cite[Chapter~22]{Snyder1987}.  Let $p$ be a point of $\sigma$ and $\theta$ the angle
between $p$ and $c$, so that $\cos\theta=|c|/|p|$.  The gnomonic projection
stretches lengths at $\chi(p)$ by $1/\cos^2\theta$ in the direction of
$p-c$ and by $1/\cos\theta$ perpendicular to it \cite[(22-2) and (22-3)]{Snyder1987}.
So $\chi$ is quasiconformal on $\sigma$, with dilatation $1/\cos\theta$ at $p$
\cite[Section~I.9.4]{LehtoVirtanen1973}, and the round metric pulled back by
$\chi$ to the plane of $\sigma$ is, at $p$, an inner product $G'$ with
$(\cos^4\theta/|c|^2)\,G\le G'\le(\cos^2\theta/|c|^2)\,G$, where $G$ is the
Euclidean inner product.  All faces are congruent, and on the
face with vertices proportional to $(0,1,\phi)$, $(0,-1,\phi)$, $(\phi,0,1)$ we have
$|c|^2=(3\phi+2)/(3\phi+6)$ and $|p|\le1$, so $\cos\theta\ge|c|>0.79$.  The faces are
equilateral, so by Lemma~\ref{lem: angle comparison} their angles in $G'$ are at most
$\pi/3+2\arcsin(0.12)<74^\circ$.  So $K_0$ is pointwise acute.

\subsubsection{Flat tori}\label{sec: acute torus}
On a flat torus the triangles can be parametrised affinely, so an acute geodesic
triangulation is pointwise acute.  The flat torus $\C/\Lambda$ has a
triangulation into $16$ acute geodesic triangles
\cite[Section~2]{ItohYuan2009}.

\subsection{Genus at least two}\label{sec: cdvm appendix}
In the final case for Theorem~\ref{prop: acute base exists}, genus at least
two, we build the base triangulation from triangles that are equilateral in a
flat metric.  Following Colin de Verdi\`ere and Marin
\cite[Section~4]{ColinDeVerdiereMarin1990}, we replace $g$ by a nearby
metric $g'$ that is conformal to a flat metric with cone points in which
$\Sigma$ can be cut into equilateral triangles
(Section~\ref{sec: cdvm conformally equilateral}).  The affine
parametrisation of these triangles in the flat coordinates is not smooth at
a cone point, and Section~\ref{sec: cdvm cone points} smoothly reparametrises the
triangles there by an explicit diffeomorphism that keeps their angles acute.

\subsubsection{Conformally equilateral metrics}\label{sec: cdvm conformally equilateral}
We use the complex structure associated with a metric and the orientation of
$\Sigma$.  The lattice $L\Z[e^{i\pi/3}]$ is the set of vertices of a grid of equilateral
triangles of side $L$ in $\C$.  Our triangulation comes from overlaying this lattice on $\Sigma$. To do so, we look for flat metrics with cone points of angle $7\pi/3$ that are conformal to a $g'$ near $g$. Such flat metrics come from sextic differentials.  Let $\Phi'$ be a sextic differential with simple zeros, holomorphic for $g'$.  Away from the zeros, the integrals of the local sixth roots of $\Phi'$ are holomorphic coordinates, unique up to $z\mapsto\xi z+\beta$ with $\xi^6=1$, and the Euclidean metric in these coordinates is a flat metric
conformal to $g'$, with a cone point of the desired angle at each zero
\cite[Section~4]{ColinDeVerdiereMarin1990}.  Fix a zero $p_0$ and $L>0$.
 The integrals along loops based at $p_0$ and along
paths from $p_0$ to the other zeros lie in $L\Z[e^{i\pi/3}]$ exactly when that grid
can be laid over $\Sigma$ in the flat coordinates with the zeros of $\Phi'$ at
vertices, and the grid then gives a triangulation of $\Sigma$
\cite[Section~4]{ColinDeVerdiereMarin1990},
\cite[(11) in the proof of Proposition~6.2]{Engel2021}, with seven triangles at each zero
filling the cone angle $7\pi/3$.  The metrics $g'$ we need are those for
which this happens.

\begin{definition}[{\cite[Section~4]{ColinDeVerdiereMarin1990}}]\label{def: cdvm conformally equilateral}
    A metric $g'$ on $\Sigma$ is \emph{conformally equilateral} if there are a
    sextic differential $\Phi'$, holomorphic for $g'$ and with simple zeros, a
    zero $p_0$ of $\Phi'$ and $L_0>0$ such that all integrals of a sixth root
    of $\Phi'$ along loops based at $p_0$ and along paths from $p_0$ to the
    other zeros lie in $L_0\Z[e^{i\pi/3}]$.
\end{definition}

Colin de Verdi\`ere and Marin prove that such metrics are dense
\cite{ColinDeVerdiereMarin1990}, adapting an argument of Tischler
\cite{Tischler1970}.  We record that their proof gives the pointwise closeness
of metrics that we need.

\begin{lemma}\label{lem: cdvm density}
	For every $\epsilon>0$ there is a conformally equilateral metric $g'$ on
	$\Sigma$ with $(1-\epsilon)g\le g'\le(1+\epsilon)g$.
\end{lemma}

\begin{proof}
	Let $J$ be the complex structure of $g$.  Colin de Verdi\`ere and Marin
	obtain conformally equilateral complex structures $J'$ on $\Sigma$ from
	arbitrarily small perturbations of a sixth root of a sextic differential,
	holomorphic for $J$, by forms that vanish near its zeros
	\cite[Section~5]{ColinDeVerdiereMarin1990}.  Away from the zeros $J'$ is
	determined pointwise by the perturbed form, and near them it equals $J$,
	so $J'$ is uniformly close to $J$ for small perturbations.

	The metric $g'(u,v):=\tfrac12\big(g(u,v)+g(J'u,J'v)\big)$ satisfies
	$g'(J'u,J'v)=g'(u,v)$, as $J'^2=-1$, so its complex structure is $J'$ and
	it is conformally equilateral.  As $g(Ju,Jv)=g(u,v)$, the difference
	$g'(u,v)-g(u,v)=\tfrac12\big(g(J'u,J'v)-g(Ju,Jv)\big)$ tends to zero
	uniformly for $g$-unit vectors $u$, $v$, so
	$(1-\epsilon)g\le g'\le(1+\epsilon)g$ once $J'$ is close enough to $J$.
\end{proof}

\subsubsection{Triangulation for the conformally equilateral metric}\label{sec: cdvm cone points}
This subsubsection constructs, for a conformally equilateral metric $g'$, the
triangulation used in the proof of
Theorem~\ref{prop: acute base exists}.  It is a fine equilateral
triangulation, parametrised affinely in the flat coordinates away from the
cone points.  At a cone point the affine parametrisation is not smooth.  We
replace it by one constructed in
Lemma~\ref{lem: cone point parametrisation} below, which is smooth, fits with those
of the neighbouring triangles and keeps the angles acute, with a margin that
does not depend on $g'$.

Let $\Phi'$, $p_0$ and $L_0$ be as in
Definition~\ref{def: cdvm conformally equilateral} for $g'$, so that the grid
of side $L_0$ gives a triangulation of $\Sigma$ with a vertex at every zero
of $\Phi'$ (Section~\ref{sec: cdvm conformally equilateral}).  Dividing $L_0$
by an integer $j\ge1$ gives another triangulation with a vertex at every
zero, which we call the triangulation of side $L:=L_0/j$.  It is clearly
simplicial for $L$ small enough.

Near a zero the flat metric has a cone point of angle $7\pi/3$, and if $z$ is the flat
coordinate centred there, $w:=z^{6/7}$ opens the cone out to a plane and is a holomorphic
coordinate of $\Sigma$.  As $z$ and $w$ are holomorphic, angles
measured in $g'$ are the Euclidean angles in these coordinates.  For $L$
small enough the seven triangles at each zero of $\Phi'$ lie in a disc around it on
which $w$ is a coordinate and which contains no other zero; in
particular no two zeros are adjacent.  We fix $L$ small enough for this and
for the triangulation of side $L$ to be simplicial.

In the coordinate $w$ a triangle at a zero is a curved triangle: its
two sides at the zero are straight, and its third side is the image of a
straight segment under $z\mapsto z^{6/7}$, slightly bent.  The neighbours
across the two straight sides are triangles at the same zero, and, as no two
zeros are adjacent, the neighbour across the third side is parametrised
affinely in $z$ along it.  So we
parametrise the curved triangle from the straight triangle with the same
vertices, by a map that is the identity on the two sides at the zero and runs
along the third side affinely in $z$.  We obtain it by bending the third side
of the straight triangle onto the curved one and spreading the bend into the
interior.  Such maps go back to the curved finite elements of Zl\'amal
\cite{Zlamal1973}, see also \cite{Lenoir1986,Bernardi1989}.  The bend is slight, so the map is close to the identity, and the
angles in $g'$ stay below $77^\circ$.  This is an explicit computation.

\begin{lemma}\label{lem: cone point parametrisation}
    Let $D$ be a triangle at a zero, which is straight in the coordinate $z$
    and curved in $w$, and let $T$ be the triangle with the same vertices that is straight in the
    coordinate $w$.  In $w$ there is a diffeomorphism $\psi$
    from $T$ onto $D$, extending to a smooth embedding of a neighbourhood of $T$,
    that is the identity on the two sides of $T$ at the zero and maps the
    third side of $T$ onto that of $D$ affinely in $z$, such that for
    every $x\in T$ every angle of $T$ in the inner product $(\psi^*e)_x$ is
    less than $77^\circ$, where $e$ is the Euclidean inner product.
\end{lemma}

\begin{proof}
    After a rotation and a scaling of $z$, and the corresponding ones of
    $w$, which preserve the claim, the triangle in $z$ is the triangle
    $E$ with vertices $0$, $1$ and $e^{i\pi/3}$.  Then $D=f(E)$, where
    $f(z):=z^{6/7}$ on $\C\setminus(-\infty,0]$, with
    $\arg f(z)=\frac67\arg z$ for $-\pi<\arg z<\pi$.  The line through $1$
    and $e^{i\pi/3}$ misses $(-\infty,0]$, so
    $\rho(t):=f\big((1-t)+te^{i\pi/3}\big)$ is a smooth curve defined for
    all real $t$.  With $q:=f(e^{i\pi/3})=e^{2\pi i/7}$, the curved triangle
    $D$ is bounded by the segments $[0,1]$ and $[0,q]$ and by the segment of
    $\rho$ from $\rho(0)=1$ to $\rho(1)=q$, and the straight triangle $T$ has vertices $0$, $1$ and $q$,
    with angles $2\pi/7$ at $0$ and $5\pi/14$ at $1$ and at $q$.  So $\psi$ is
    to be the identity on $[0,1]$ and $[0,q]$ and to satisfy
    $\psi\big((1-t)+tq\big)=\rho(t)$ for $t\in[0,1]$.

    Let $c(t):=(1-t)+tq$ parametrise the third side of $T$, so that the bend
    $u:=\rho-c$ is smooth and vanishes at $t=0$ and $t=1$.  Put
    \[
        \tilde u(t):=-\iint_{a,b\ge0,\ a+b\le1}u''(a+bt)\,da\,db,
    \]
    which is smooth in $t$ and, as one checks by integrating,
    equals $u(t)/\big(t(1-t)\big)$ for $0<t<1$.  Let $\lambda_0$,
    $\lambda_1$, $\lambda_2$ be the barycentric coordinates with respect to
    $0$, $1$, $q$, let $s:=\lambda_2+\lambda_0/2$, and put
    \[
        \psi(x):=x+\lambda_1(x)\lambda_2(x)\,\tilde u\big(s(x)\big).
    \]
    On $[0,1]$ and on $[0,q]$ one of $\lambda_1$, $\lambda_2$ vanishes, so
    $\psi$ is the identity there.  At $c(t)$ we have $\lambda_1=1-t$ and
    $\lambda_2=s=t$, so $\psi(c(t))=c(t)+u(t)=\rho(t)$.

    The map $\psi$ is close to the identity because the bend is slight.  Its
    derivative is
    \[
        d\psi-\mathrm{id}
        =\tilde u(s)\,(\lambda_2\,d\lambda_1+\lambda_1\,d\lambda_2)
        +\lambda_1\lambda_2\,\tilde u'(s)\,ds .
    \]
    We bound the factors that depend only on the shape of $T$, then
    $\tilde u$ and $\tilde u'$, and combine; as $s$ takes values in $[0,1]$
    on $T$, bounds for $\tilde u$ and $\tilde u'$ on $[0,1]$ suffice.

    On $T$ we have $\lambda_1+\lambda_2\le1$ and $\lambda_1\lambda_2\le\frac14$.
    One checks that $|d\lambda_1|=|d\lambda_2|=1/\sin(2\pi/7)<4/3$, the
    reciprocals of the heights of $T$ at $1$ and at $q$, and that
    $|ds|=1/\big(2\sin(\pi/7)\big)<6/5$.

    The bounds on $\tilde u$ and $\tilde u'$ come from bounds on $u''$ and
    $u'''$.  Put $\varsigma(t):=(1-t)+te^{i\pi/3}$, so that
    $\rho=f\circ\varsigma$.  As $c$ and $\varsigma$ are affine,
    $u''=f''(\varsigma)(e^{i\pi/3}-1)^2$ and
    $u'''=f'''(\varsigma)(e^{i\pi/3}-1)^3$, and as $|e^{i\pi/3}-1|=1$ and
    $|\varsigma|\ge\sqrt3/2$, one checks that $|u''|<0.15$ and $|u'''|<0.2$ on
    $[0,1]$.  Differentiating the formula for $\tilde u$ under the integral
    brings out a factor $b$, and over the triangle of integration
    $\iint da\,db=\frac12$ and $\iint b\,da\,db=\frac16$, so on $[0,1]$
    \[
        |\tilde u|<0.15\cdot\tfrac12=0.075
        \quad\text{and}\quad
        |\tilde u'|<0.2\cdot\tfrac16=\tfrac1{30}.
    \]

    Combining these, with $\lambda_1+\lambda_2\le1$ in the first term and
    $\lambda_1\lambda_2\le\frac14$ in the second, gives on $T$
    \[
        |d\psi-\mathrm{id}|<0.075\cdot\tfrac43+\tfrac14\cdot\tfrac1{30}\cdot\tfrac65=0.11 .
    \]

    Consequently $\psi$ is a local diffeomorphism, and it is injective on the
    convex set $T$, since $|\psi(x)-\psi(y)|\ge0.89|x-y|$ by the mean value
    inequality.  By continuity $|d\psi-\mathrm{id}|<1$ on a convex open
    neighbourhood of $T$, and the same argument makes $\psi$ injective there, so
    $\psi$ is a smooth embedding of that neighbourhood.  As $\psi$ maps the boundary of
    $T$ onto that of $D$, it maps $T$ onto $D$.  Finally $0.89^2e\le(\psi^*e)_x\le1.11^2e$, so by
    Lemma~\ref{lem: angle comparison} the angles of $T$ in $(\psi^*e)_x$
    exceed its Euclidean angles, at most $5\pi/14$, by at most
    $2\arcsin(0.11)$, and $5\pi/14+2\arcsin(0.11)<77^\circ$.
\end{proof}

We now parametrise the triangles of the triangulation of side $L$: affinely
in $z$ away from the zeros, with angles $\pi/3$ in $g'$, and at a zero by
the affine map onto the straight triangle with the same vertices in
$w$, followed by a map $\psi$ as in
Lemma~\ref{lem: cone point parametrisation}, with angles less than
$77^\circ$ in $g'$.  These
parametrisations agree on shared edges, by the boundary conditions in that
lemma, so they form a smooth pointwise acute triangulation $(K_0,\chi)$ of $(\Sigma,g')$, with a margin that does not depend on $g'$.

By Lemma~\ref{lem: cdvm density}, for any $\epsilon>0$, we can choose the conformally equilateral metric $g'$ on $\Sigma$ such that $(1-\epsilon)g\le g'\le(1+\epsilon)g$.     Passing from $g'$ back to $g$ changes the angles of $(K_0,\chi)$ only a little: they change by at most the bound~\eqref{eq: angle distortion} with $a^2=1-\epsilon$ and $b^2=1+\epsilon$ (Lemma~\ref{lem: angle comparison}).  So for $\epsilon$ small enough $(K_0,\chi)$ is a smooth pointwise acute triangulation of $(\Sigma,g)$.

\section*{Acknowledgements}
The author is grateful to Owen Gwilliam for encouragement while preparing this manuscript.
This research was supported by the European Union's Horizon Europe
research and innovation programme under the Marie Sk\l{}odowska-Curie
grant agreement No~101203556 (DisQ FA), and by Aalto University.

\section*{Tool and computational resource disclosure}
Large language models (predominantly Claude) were used throughout the drafting of this paper. They were used for literature search and for drafting and revising prose and proofs. The author has revised and edited all text, and has ensured that every proof is correct, readable and properly attributed.

The author takes full responsibility for the correctness of the arguments and results, and for the completeness and accuracy of the citations.  No automated system is an author of this paper; the core intellectual contributions are the author's.

\bibliographystyle{alpha}
\bibliography{lib}

\end{document}